\documentclass[12pt]{amsart}
\usepackage{amssymb, amsmath, amscd, mathrsfs, marvosym, psfrag} 

\input xy
\xyoption{all}
\DeclareMathAlphabet{\mathpzc}{OT1}{pzc}{m}{it}

\newtheorem{theorem}{Theorem}[section]
\newtheorem{proposition}[theorem]{Proposition}
\newtheorem{corollary}[theorem]{Corollary}

\newtheorem{conjecture}[theorem]{Conjecture}
\newtheorem{lemma}[theorem]{Lemma}
\newtheorem{definition}[theorem]{Definition}

\theoremstyle{remark}
\newtheorem{remark}[theorem]{Remark}

\newtheorem{Problem}[theorem]{Problem}

\def\varle{\leqslant}
\def\varge{\geqslant}

\newcommand{\CA}{{\mathcal A}}

\newcommand{\CC}{{\mathcal C}}

\newcommand{\CE}{{\mathcal E}}
\newcommand{\CF}{{\mathcal F}}
\newcommand{\CG}{{\mathcal G}}

\newcommand{\CI}{{\mathcal I}}
\newcommand{\CJ}{{\mathcal J}}

\newcommand{\CO}{{\mathcal O}}
\newcommand{\CP}{{\mathcal P}}

\newcommand{\CS}{{\mathcal S}}

\newcommand{\CV}{{\mathcal V}}
\newcommand{\CW}{{\mathcal W}}

\newcommand{\CZ}{{\mathcal Z}}

\newcommand{\SB}{{\mathscr B}}

\newcommand{\SF}{{\mathscr F}}
\newcommand{\SG}{{\mathscr G}}

\newcommand{\SV}{{\mathscr V}}

\newcommand{\SZ}{{\mathscr Z}}

\newcommand{\fh}{{{\mathfrak h}}}

\newcommand{\fg}{{{\mathfrak g}}} 
\newcommand{\fb}{{{\mathfrak b}}}

\newcommand{\fm}{{{\mathfrak m}}}

\newcommand{\fhd}{\fh^\star}

\newcommand{\hfhd}{\widehat{\fh}^\star}

\newcommand{\tS}{\widetilde{S}}

\newcommand{\DC}{{\mathbb C}}
\newcommand{\DP}{{\mathbb P}}
\newcommand{\DR}{{\mathbb R}}
\newcommand{\DZ}{{\mathbb Z}}

\newcommand{\DV}{{\mathbb V}}

\newcommand{\DW}{{\mathbb W}}

\newcommand{\ch}{{\operatorname{ch}\, }}

\newcommand{\End}{{\operatorname{End}}}

\newcommand{\Hom}{{\operatorname{Hom}}}

\newcommand{\id}{{\operatorname{id}}}

\newcommand{\supp}{{\operatorname{supp}}}
\newcommand{\catmod}{{\operatorname{-mod}}}

\newcommand{\catmodr}{{\operatorname{mod-}}}

\newcommand{\rk}{{{\operatorname{rk}}}}

\newcommand{\ol}{\overline}

\newcommand{\pt}{{pt}}

\newcommand{\linie}{{\,\text{---\!\!\!---}\,}}
\newcommand{\llinie}{{\text{---\!\!\!---\!\!\!---}}}

\newcommand{\GL}{{\operatorname{GL}}}

\newcommand{\IC}{{\operatorname{IC}}}

\newcommand{\ul}{\underline}

\newcommand{\THypgr}{{\operatorname{\mathbb{H}_T^\bullet}}}

\newcommand{\HHypgr}{{\operatorname{\mathbb{H}_H^\bullet}}}

\newcommand{\comment}[1]{}
\begin{document}

\pagenumbering{arabic}
\title[Moment graphs]{Moment graphs in representation theory and geometry}
\author[]{Peter Fiebig, Universit\"at Erlangen-N\"urnberg}
\begin{abstract} This paper reviews the moment graph technique that allows to translate certain representation theoretic problems into geometric ones. For simplicity we restrict ourselves to the case of semisimple complex Lie algebras. In particular, we show how the original Kazhdan--Lusztig conjecture on the characters of irreducible highest weight representations can be translated into a multiplicity problem for parity sheaves on the (Langlands dual) flag variety. 

\end{abstract}
\maketitle
\section{Introduction}
One of the central problems in representation theory is the determination of the simple characters of a given group or, more generally, of a given algebra.  Often this problem turns out to be difficult to solve and there is an abundance of situations in which we do not  have a sufficient answer. The problem seems to be particularly difficult if the base field of our theory  is of prime characteristic. 

A very successful approach towards the character problem is to find a relation between representation theory and the geometry of some algebraic variety. If such a relation is established one might hope that the machinery of algebraic geometry is powerful enough to provide a solution.

There are at least two  different ways to relate representation theory to geometry. The first, which  well-known and well established in several areas, appeared for the first time in the article \cite{BB} and is called the   {\em Beilinson--Bernstein localization}. It amounts to realizing a Lie algebra inside the space of differential operators on a complex algebraic variety. It thus allows us  to relate the category of $D$-modules to representations of the Lie algebra via the global sections functor. A characteristic $p$ analogue of this approach has been studied by Bezrukavnikov et al. in \cite{BMR}.

The second idea goes back to Soergel. Here,  the geometry and the representation theory are related via an intermediate category, which often is of a linear-algebraic, ``combinatorial'' nature. Geometrically it arises via the hypercohomology functor applied to an appropriate class of sheaves (such as intersection cohomology complexes or parity sheaves). Representation theoretically this category is the image of versions of Soergel's {\em Strukturfunktor}. 

These ``combinatorial'' categories can sometimes be realized over arbitrary fields and even over the integers, which opens a way to compare the structure over different fields using base change results (see, for example, \cite{AJS}). Here we want to discuss the main ideas of Soergel's approach in the particular example of representations of a semisimple complex Lie algebra, which is technically the least demanding instance. The associated combinatorial category comes in two quite different and important realizations: the first is the category of Soergel bimodules, the second is the category of sheaves on a moment graph. In the complex semisimple case this is a finite, directed and labelled graph that can be read of from the underlying root system and its Weyl group. 

Following Soergel's idea, we want to state the main steps and explain the main ideas for a proof of the classical Kazhdan-Lusztig conjecture using moment graphs. A characteristic $p$ version of this approach can be found in the articles \cite{FieModRep, FieCharIrr}.

While both approaches explained  above typically yield character formulas in various situations,  they also relate categorical structures.  Their real power is revealed once one combines the two: this  establishes the celebrated {\em Koszul self duality} of category $\CO$.

The first section of this article is the most elementary one and deals with the basics of moment graph theory. In particular, we introduce the moment graph associated to a finite root system. Then we motivate the construction of the principal objects associated to a moment graph: the Braden--MacPherson sheaves.

The second section  describes the link between the representation theory of a finite dimensional complex Lie algebra and the moment graph theory. We present the most classical case of the Kazhdan--Lusztig conjecture and show how it can be translated into a multiplicity conjecture on the stalks of the Braden--MacPherson sheaves introduced earlier. 

The third and final section gives the link between the geometry of flag varieties and moment graph theory. We present the localization functor that associates a sheaf on the moment graph to a torus-equivariant sheaf on the variety following \cite{GKM}. Here we also cover the positive characteristic case following \cite{FW}. The most important class of equivariant sheaves are the {\em parity sheaves} introduced by Juteau, Mautner and Williamson in \cite{JMW}. By localizing parity sheaves one obtains the Braden--MacPherson sheaves on the graph. Deep results in algebraic geometry show that the parity sheaves that we consider are {\em intersection cohomology sheaves} if the characteristic of the field of coefficients is either zero or large enough\footnote{A direct link between intersection cohomology sheaves with coefficients in $\DC$ and Braden--MacPherson sheaves is established in \cite{BMP}.}. A formula for the dimension of the stalks of the latter sheaves then yields  multiplicities for Braden--MacPherson sheaves and, in particular, proves the Kazhdan--Lusztig conjecture we presented in the second section.

\section{Moment graphs}

We  introduce the notion of a moment graph and then present the most important class of moment graphs, the moment graphs associated to a root system. We define sheaves on such graphs and the corresponding functor of local sections. A natural extension  property then leads  to the definition of the Braden--MacPherson sheaves. We also  discuss another characterization of the Braden--MacPherson sheaves: their global sections form a set of representatives of the projective  objects in a certain exact category. This is the property that latter allows us to give these sheaves a representation theoretic meaning. A reference for most of the following is \cite{FieAdv}.

\subsection{Moment graphs}\label{subsec-momgra}
Let $Y$ be a lattice of finite rank (i.e.~$Y\cong\DZ^r$).
An (unordered) moment graph over $Y$ is the datum $\CG=(\CV,\CE,l)$, where $(\CV,\CE)$ is a finite graph without loops (edges connecting a vertex to itself) and no double edges (two vertices are connected by at most one edge), and $l\colon \CE\to Y\setminus\{0\}$ is a map. We write $E\colon x\linie y$ to denote an edge connecting $x$ and $y$ and we write $E\colon x\stackrel{l(E)}\llinie y$ if we also want to specify its label. 

One obtains an important class of moment graphs from root systems in the following way. Let $R\subset V$ be a root system in a real vector space $V$  and let us, for latter use, fix a system of positive roots $R^+\subset R$. For each $\alpha\in R$ we denote by $\alpha^\vee\in V^\ast=\Hom_\DR(V,\DR)$ its coroot. By $R^\vee\subset V^\ast$ we denote the coroot system. By $Y=\{v\in V^\ast\mid \alpha(v)\in\DZ\text{ for all $\alpha\in R$}\}$ we denote the coweight lattice. By $\CW\subset\GL(V)$ we denote the Weyl group that is generated by the reflecions $s_\alpha$ associated to the roots $\alpha\in R^+$.

We now define the graph $\CG_R$ over the coweight lattice $Y$:
\begin{itemize}
\item its set of vertices is $\CW$,
\item the elements $x,y\in\CW$ are connected by an edge if there is a root $\alpha\in R^+$ with $x=s_{\alpha} y$,
\item for $E\colon x\linie s_{\alpha} x$ and $\alpha\in R^+$ we set $l(E)=\alpha^\vee$.
\end{itemize}

With any moment graph $\CG$ and any field $k$ we can now associate the category of {\em sheaves with coefficients in $k$}.

\subsection{Sheaves on moment graphs} Suppose that $\CG$ is a moment graph defined over the lattice $Y$. 
Now let us fix a field $k$. We denote by $Y_k=Y\otimes_\DZ k$ the $k$-vector space spanned by  $Y$. By $S=S(Y_k)$ we denote its symmetric algebra. This is a polynomial algebra (over $k$) of rank $\dim_k Y_k$ and we consider it as a $\DZ$-graded algebra with $Y_k\subset S$ being the homogeneous component of degree $2$. By an $S$-module we mean a graded $S$-module in the following. The shift functor $[1]$ on a graded $S$-module $M$ is defined in such a way that $M[1]_l=M_{l+1}$  for each $l\in\DZ$.

A {\em $k$-sheaf} on $\CG$ is the datum $\SF=(\SF^x,\SF^E,\rho_{x,E})$, where
\begin{itemize}
\item $\SF^x$ is an $S$-module for any vertex $x\in\CV$,
\item $\SF^E$ is an $S$-module with $l(E)\SF^E=0$ for each edge $E\inÊ\CE$,
\item $\rho_{x,E}\colon \SF^x\to \SF^E$ is a homomorphism of $S$-modules for any vertex $x$ adjacent to the edge $E$.
\end{itemize}
A morphism $f\colon \SF\to\SG$ of sheaves on $\CG$ is given by homomorphisms of $S$-modules $f^x\colon \SF^x\to\SG^x$ and $f^E\colon \SF^E\to \SG^E$ for all $x$ and $E$ such that for any vertex $x$ lying on the edge $E$ the diagram

\centerline{
\xymatrix{
\SF^x\ar[d]_{\rho_{x,E}}\ar[r]^{f^x}&\SG^x\ar[d]{\rho_{x,E}}\\
\SF^E\ar[r]^{f^E}&\SG^E
}
}
\noindent
commutes. We denote by $\CG\catmod_k$ the corresponding category of $k$-sheaves on $\CG$. This category is graded by the shift functor $[1]\colon \CG\catmod_k\to\CG\catmod_k$ that applies the functor $[1]$ to each component of a sheaf.

Here are some examples of sheaves on moment graphs.
\begin{enumerate}
\item The most natural sheaf on $\CG$ is probably the {\em structure sheaf} $\SZ$ that is defined by $\SZ^x=S$ for all $x\in\CV$ and $\SZ^E=S/l(E) S$ for any edge $E$, and $\rho_{x,E}\colon S\to S/l(E) S$ the canonical quotient map.
\item For each vertex $x\in\CV$ there is a {\em skyscraper sheaf} $\SV(x)$ defined by $\SV(x)^x=S$, and all other data is zero.
\end{enumerate}

\subsection{Sections of sheaves}

Let $\SF$ be a sheaf on the moment graph $\CG$ and let $\CI$ be a a subset of $\CV$. The space of sections of $\SF$ over $\CI$ is
$$
\Gamma(\CI,\SF)=\left\{
(f_x)\in\bigoplus_{x\in\CI} \SF^x\left|
\,
\begin{matrix}
\rho_{x,E}(f_x)=\rho_{y,E}(f_y) \\
\text{for all edges $E\colon x\linie y$ with $x,y\in\CI$}
\end{matrix}
\right\}.
\right.
$$
In particular, we define $\Gamma(\SF):=\Gamma(\CV,\SF)$ as the space of global sections. For each pair of subsets $\CI$,  $\CI^\prime$ of $\CV$ with $\CI^\prime\subset \CI$ the projection $\bigoplus_{x\in\CI}\SF^x\to\bigoplus_{x\in\CI^\prime}\SF^x$ along the decomposition induces a restriction map $\Gamma(\CI,\SF)\to\Gamma(\CI^\prime,\SF)$.

The space 
$$
\CZ=\Gamma(\SZ)=
\left\{
(f_x)\in\bigoplus_{x\in\CV} S
\left|Ê
\begin{matrix} 
f_x\equiv f_y\mod l(E) \\
\text{ for all edges $E\colon x\linie y$}
\end{matrix}
\right.
\right\}
$$
of global sections of the structure sheaf carries a canonical structure of a graded $S$-algebra. The addition and the multiplication are given componentwise. The space of global sections of any sheaf $\SF$ naturally forms a graded $\CZ$-module.
\begin{remark} It is possible to consider a moment graph sheaf as we defined it above as a sheaf on a topological space in the usual sense (see \cite{BMP}). However, for the following this viewpoint seems not to be helpful yet. 
\end{remark}
\subsection{A topology on the moment graph}
A natural problem in sheaf theory is the question whether a sheaf is {\em flabby}, i.e.~whether the restriction map from the global sections to the sections over an open subset is surjective. We want to study the analogous question in our version of  sheaf theory, so we need a notion of an open subset. 

For this we now assume that we are given an additional structure on the moment graph $\CG$, namely we assume that each edge of $\CG$ is directed. We write $E\colon x\to y$ or $E\colon x\stackrel{\alpha}\to y$ for a directed edge.  For $x,y\in \CV$ we then set $x\varle y$ if either $x=y$ or if there is a directed path leading from $x$ to $y$, i.e.~ if there are vertices $x_0,\dots, x_n$ with $x=x_0$, $y=x_n$ and directed edges $x_{i-1}\to x_{i}$ for all $i=1,\dots, n$. We assume that ``$\varle$'' defines a partial order on $\CV$, i.e.~ we assume that there are no directed cycles. 

In the case of the moment graph associated to the root system $R$ we use the following directions: For an edge $E\colon x\to y$ there is a unique positive root $\alpha$ with $x=s_\alpha y$. We direct $E$ towards $y$ if $y$ is bigger than $x$ in the Bruhat order.

Let us return to the case of a general moment graph. 
The partial order  can be used to define a topology on the set $\CV$: a subset 
$\CI$ of $\CV$ is {\em open} if for any $x\in\CI$ and  $y\in\CV$ with $x\varle y$ we have $y\in\CI$, i.e.~ if $\CI=\bigcup_{x\in\CI}\{\varge x\}$. Here and in the following we write $\{\ge z\}$ for $\{y\in \CV\mid y\ge z\}$. The notation $\{>z\}$ has an analogous meaning.
\subsection{A restriction on the characteristic}
An essential assumption for us is that the characteristic of $k$ is not too small with respect to the labels on the graph:
\begin{definition} Let $\CG$ be a moment graph and $k$ a field. We say that $(\CG,k)$ {\em satisfies the GKM-assumption} if $\ch k\ne 2$ and if for any distinct edges $E$ and $E^\prime$ adjacent at a common vertex we have $l(E)\not\in kl(E^\prime)$ (this is to be read in $Y\otimes_\DZ k$).
\end{definition}

In case the moment graph is associated to the (finite)  root system $R$, we only have to exclude characteristic $2$ and characteristic $3$ in case $R$ has a component of type $G_2$. For subgraphs of affine moment graphs, however, the above is a crucial restriction (see \cite{FieModRep}).

\subsection{The Braden-MacPherson sheaf}
Now we want to ask whether a given sheaf $\SF$ is {\em flabby}, i.e.~ whether for any open subset $\CI$ of $\CV$ the restriction of the global sections 
$\Gamma(\SF)\to \Gamma(\CI,\SF)$ is surjective. There is a certain class of sheaves (cf.~\cite{BMP}) that are universal with respect to the problem of extension of local sections:

\begin{definition}[\cite{BMP, FieAdv, FW}] A sheaf $\SB$ on the moment graph $\CG$ is called a {\em Braden--MacPherson sheaf} if it satisfies the following properties:
\begin{enumerate}
\item $\SB^x$ is a graded free $S$-module of finite rank for any $x\in\CV$,
\item for a directed edge $E\colon x\to y$ the map $\rho_{y,E}\colon \SB^y\to\SB^E$ is surjective with kernel $l(E)\SB^y$,
\item for any open subset $\CJ$ of $\CV$ the map $\Gamma(\SB)\to \Gamma(\CJ,\SB)$ is surjective, and
\item the composition $\Gamma(\SB)\subset\bigoplus_{z\in\CV}\SB^z\to \SB^x$ is surjective for any $x\in \CV$ (where the  map on the right is the projection along the decomposition).
\end{enumerate}
\end{definition}

The property (4) is a technical property that assures that the stalks are as small as they need to be. One can quite easily deduce the following results on the category of Braden--MacPherson sheaves.

\begin{theorem} [\cite{BMP,FieAdv,FW}]
\begin{enumerate}
\item For any $w\in\CV$ there is an up to isomorphism unique Braden--MacPherson sheaf $\SB(w)$ on $\CG$ with the following properties:
\begin{itemize}
\item We have $\SB(w)^w\cong S$ and $\SB(w)^x= 0$ unless $x\le w$.
\item $\SB(w)$ is indecomposable in $\CG\catmod_k$.
\end{itemize}
\item Let $\SB$ be a Braden--MacPherson sheaf. Then there are $w_1,\dots,w_n\in \CV$ and $l_1,\dots,l_n\in\DZ$ such that
$$
\SB\cong\SB(w_1)[l_1]\oplus \dots\oplus\SB(w_n)[l_n].
$$
The multiset $(w_1,l_1)$,\dots,$(w_n,l_n)$ is uniquely determined by $\SB$.
\end{enumerate}
\end{theorem}

For latter use we now study the connection between the category of sheaves and the category of $\CZ$-modules obtained as global sections of sheaves.

\subsection{A generic decomposition of $\CZ$-modules}

Again we return to the case of a general moment graph.
Let $Q$ be the quotient field of $S$. The embedding $\CZ\subset\bigoplus_{x\in\CV}S$ induces an embedding
$
\CZ\otimes_S Q\subset \bigoplus_{x\in\CV}Q.
$
The following is quite easy to prove in the case that the moment graph is finite.
\begin{lemma}[{\cite[Lemma 3.2]{FieAdv}}] The latter embedding is a bijection, i.e.~ we have
$
\CZ\otimes_S Q=\bigoplus_{x\in\CV} Q.
$
\end{lemma}

Let $M$ be a $\CZ$-module. We denote by $M_Q$ the $Q$-vector space $M\otimes_S Q$ obtained from $M$ by base change. It is naturally acted upon by $\CZ\otimes_S Q=\bigoplus_{x\in\CW} Q$, hence there is a natural decomposition
$
M_Q=\bigoplus_{x\in\CW} M_Q^x,
$
such that $(z_x)\in\bigoplus_{x\in\CW} Q$ acts on $M_Q^x$ as multiplication with $z_x$. 

\begin{definition}[{\cite[Definition 2.7]{FieTAMS}}] Let $\CI\subset \CV$ be a subset and $M$ a $\CZ$-module. We let $M^\CI$ be the image of the composition $M\to M_Q\to \bigoplus_{x\in\CI} M_Q^x$, where the map on the left is the canonical map that sends $m\in M$ to $m\otimes 1$ in $M_Q$, and the map on the right is the projection with kernel $\bigoplus_{x\in\CV\setminus\CI} M_Q^x$.
\end{definition} 
For a vertex $x$ we write $M^x$ instead of $M^{\{x\}}$ and call this space the {\em stalk of $M$ at $x$}.

The following statement readily follows from part (4) of the definition of a Braden--MacPherson sheaf. It allows us to determine the stalks of the Braden--MacPherson sheaves from the spaces of their global sections.
\begin{lemma} Let $\SB$ be a Braden--MacPherson sheaf. For each $x\in\CV$ we have a canonical identification
$
\Gamma(\SB)^x\cong \SB^x.
$
\end{lemma}

We now put the sheaf theory in a more categorical context.
\subsection{$\CZ$-modules admitting a Verma flag}
The following is a very crucial definition for us. 
\begin{definition}[{\cite[Section 4]{FieAdv}}] Let $M$ be a $\CZ$-module. We say that $M$ {\em admits a Verma flag} if for any open set $\CJ$ the module $M^\CJ$ is graded free over $S$ of finite rank. 
\end{definition}
We denote by $\CC^{Verma}\subset \CZ\catmod$ the full category that consists of all modules admitting a Verma flag. The category $\CC^{Verma}$ is not an abelian category. In the following we define an {\em exact structure} on $\CC^{Verma}$. The notion of an exact structure on an additive category $\CA$ was introduced by Quillen in \cite{Q} (it is  a collection of sequences $A\to B\to C$ in $\CA$ that satisfies certain axioms). We need the following exact structure:

\begin{definition}[{\cite[Section 4.1]{FieAdv}}] Let $A\to B\to C$ be a sequence in $\CC^{Verma}$. We say that it is {\em exact} if for any open subset $\CJ$ of $\CV$ the induced sequence
$
0\to A^\CJ\to B^\CJ\to C^\CJ\to 0
$
 is an exact sequence of abelian groups.
\end{definition}

In an exact category $\CA$ one has the usual notion of projective objects: An object $P$ in $\CA$ is projective, if the functor $\Hom_\CA(P,\cdot)\colon \CA\to \DZ\catmod$ maps the given class of short exact sequences to short exact sequences of abelian groups.
We now obtain a categorical characterization of the Braden--MacPherson sheaves:

\begin{theorem}[\cite{FieAdv,FieTAMS}] Suppose that $\CG$ is the moment graph associated to the root system $R$ and $k$ is such that $(\CG,k)$ satisfies the GKM-assumption. Then the following holds:
\begin{enumerate}
\item If $\SB$ is a Braden--MacPherson sheaf on $\CG$, then $\Gamma(\SB)$ admits a Verma flag.
\item $\Gamma$ induces an equivalence between the full subcategory of Braden--MacPherson sheaves in $\CG\catmod_k$ and the full subcategory of projective objects in $\CC^{Verma}$. 
\end{enumerate}
 \end{theorem}
 
 In particular, the set $\{\Gamma(\SB(w))\}_{w\in\CW}$ is a full set of representatives for the indecomposable projective isomorphism classes in $\CC^{Verma}$. We can now formulate the main problem in moment graph theory.
 
 \begin{Problem} Determine the graded rank of $\SB(w)^x$ for all $w,x\in\CW$.
 \end{Problem}
 
 This rank will depend on the characteristic of the chosen field $k$. We will see in the following that in the case $k=\DC$ this rank calculates the characters of simple highest weight modules for a semi-simple complex Lie algebra and that it can be determined via the geometry of flag varieties. A direct computation of this rank  (in the equivalent language  of Soergel bimodules) avoiding the passage to geometry  can be found in \cite{EW}.
 
\section{Representation theory}

The aim of this section is to discuss the original Kazhdan--Lusztig conjecture which is a prototype for a series of similar conjectures in representation theory. We start with some basic notions in the representation theory of semi-simple complex Lie algebras. A reference for a large part of the following is the book \cite{Hum}.

\subsection{Simple highest weight modules} 
Let $\fg$ be a complex simple Lie algebra and $\fh\subset\fb\subset\fg$ a Cartan and a Borel subalgebra of $\fg$. We assume that we have identified the root system associated to $\fh$ and $\fg$ with $R$ in such a way that $R^+$ is the set of roots of $\fb$. Recall that a $\fg$-module $M$ is a {\em weight module} if it is semi-simple for the action of $\fh$, i.e.~ if $M=\bigoplus_{\lambda\in\fhd} M_\lambda$, where $M_\lambda$ is the $\fh$-weight space with eigenvalue $\lambda$.

There is  a partial order ``$\varle$'' on the set $\fhd$ given by $\lambda\varle\mu$ if and only if $\mu-\lambda$ is a sum of positive roots. Let $\lambda\in\fhd$. A module $M$ is called a  {\em module of highest weight $\lambda$}  if there is $m\in M_\lambda$, $m\ne 0$,  that generates $M$ and such that $M_\mu\ne 0$ implies $\mu\varle\lambda$. 

It is not difficult to see that for any $\lambda\in\fhd$ there is an up to isomorphism unique simple module $L(\lambda)$ of highest weight $\lambda$ and that $L(\lambda)\cong L(\mu)$ implies $\lambda=\mu$.

\subsection{Characters of highest weight modules}

Let $\DZ[\fhd]$ be the lattice spanned by the set $\fhd$. We denote the basis element associated to $\lambda$ by $e^{\lambda}$. An element in $\DZ[\fhd]$ is hence a finite $\DZ$-linear combination $\sum_{\lambda\in\fhd} a_\lambda e^\lambda$. As the modules $L(\lambda)$ are infinite dimensional in general, we need the following completion of the lattice $\DZ[\fhd]$: we let $\widehat{\DZ[\fhd]}$ be the abelian group of formal $\DZ$-linear combinations $\sum_{\lambda\in\fhd}a_\lambda e^\lambda$ such that there exist $\mu_1,\dots,\mu_n$ in $\fhd$ such that $a_\lambda\ne 0$ implies $\lambda\le \mu_i$ for some $i$. 

Now, a highest weight module $M$ has finite dimensional weight spaces $M_\lambda$, hence we can define its character as
$$
\ch M:=\sum_{\lambda\in\fhd} \dim_\DC M_\lambda\cdot e^{\lambda}\in\widehat {\DZ[\fhd]}.
$$
Our principal aim is to calculate the {\em simple highest weight characters} of $\fg$, i.e.~ to give a formula for $\ch L(\lambda)$ for any $\lambda$. 

\subsection{Simple highest weight characters and the Kazhdan--Lusztig conjecture}
Set $
\rho=1/2\sum_{\alpha\in R^+}Ê\alpha
$ and define the {\em dot-operation} of $\CW$ on $\fhd$ by shifting the linear action by $-\rho$, i.e.~ such that for $w\in\CW$ and $\lambda\in\fhd$ we have
$
w.\lambda=w(\lambda+\rho)-\rho.
$

If $\lambda$ is integral (i.e.~ $\langle\lambda,\alpha^\vee\rangle\in\DZ$ for all $\alpha\in R$), regular (i.e.~ the dot-stabilizer of $\lambda$ in $\CW$ is trivial) and dominant (i.e.~ maximal in its $\CW$-dot-orbit), then $\ch L(\lambda)$ is given by {\em Weyl's character formula}:
$$
\ch L(\lambda)=\frac{\sum_{w\in\CW} (-1)^{l(w)} e^{w(\lambda+\rho)}}{\sum_{w\in\CW} (-1)^{l(w)} e^{w(\rho)}}.
$$
Kazhdan and Lusztig conjectured in \cite{KL1} the following generalization of Weyl's character formula.

\begin{conjecture}[\cite{KL1}] \label{conj-KL} Suppose that $\lambda$ is integral, regular and dominant. Then we have for  any $w\in\CW$ 
$$
\ch L(w. \lambda)= \sum_{x\in\CW} (-1)^{l(w)-l(x)} P_{x,w}(1) e^{x. \lambda}\prod_{\alpha\in R^+}(1+e^{-\alpha}+e^{-2\alpha}+\dots).
$$
Here  $P_{x,y}$ denotes the Kazhdan--Lusztig polynomial for the Coxeter system $(\CW,\CS)$ at the parameters $x,y$. 
\end{conjecture}

Our next task is to rewrite the Kazhdan--Lusztig conjecture as a multiplicity conjecture for Verma modules.

\subsection{Verma modules}

Recall that  we have $\fb=\fh\oplus[\fb,\fb]$ as vector spaces.  For any linear form $\lambda\in\fhd$ we let $\DC_\lambda$ be the $1$-dimensional $\fb$-module on which $\fh\subset \fb$ acts via the character $\lambda$ and $[\fb,\fb]$ acts trivially. Induction then yields the {\em Verma module}
$$
\Delta(\lambda)=U(\fg)\otimes_{U(\fb)} \DC_\lambda.
$$
This is a module of highest weight $\lambda$ and it even is {\em universal} in the sense that any module of highest weight $\lambda$ (in particular, $L(\lambda)$) is isomorphic to a quotient of $\Delta(\lambda)$. 
From the PBW-theorem one deduces  immediately the following character formula for the Verma modules.

\begin{lemma}\label{lemma-charVerma}  For any $\lambda\in\fhd$ we have
$$
\ch \Delta(\lambda)=e^{\lambda}\prod_{\alpha\in R^+}(1+e^{-\alpha}+e^{-2\alpha}+\dots).
$$
\end{lemma}

\subsection{Jordan--H\"older multiplicities}
Each highest weight module $M$ admits a Jordan-H\"older series, i.e.~ a filtration
$
0=M_0\subset M_1\subset \dots\subset M_n=M
$
such that $M_{i}/M_{i-1}\cong L(\lambda_i)$ for some $\lambda_1,\dots, \lambda_n\in\fhd$. For each $\mu$ the number of occurences of $L(\mu)$ in such a filtration is independent of the choice of the filtration and is denoted by $[M:L(\mu)]$.

Using Lemma \ref{lemma-charVerma} and an inversion formula for the Kazhdan--Lusztig polynomials  one shows that the following conjecture is equivalent to the Kazhdan--Lusztig conjecture \ref{conj-KL}. 
\begin{conjecture}[\cite{KL1}]\label{conj-KL2} Suppose that $\lambda$ is integral,  regular and dominant. Then we have for all $w\in\CW$
$$
[\Delta(w.\lambda):L(\mu)]=
\begin{cases}
P_{w_0w,w_0x}(1),&\text{ if $\mu=x.\lambda$ for some $x\in\CW$,}\\
0,& \text{ otherwise.}
\end{cases}
$$
\end{conjecture}
Here $w_0\in\CW$ denotes the longest element.

\subsection{The category $\CO$}
In order to study the above multiplicity conjecture  we need to define a category around our objects. This category is the ``highest weight''  category $\CO$ defined originally by Bernstein, Gelfand and Gelfand. Here is its definition:

\begin{definition} We denote by $\CO$ the full subcategory of the category of $\fg$-modules that consists of locally $\fb$-finite weight modules.
\end{definition}
Recall that a $\fg$-module $M$ is called {\em locally $\fb$-finite}, if every $m\inÊM$ is contained in a finite dimensional $\fb$-submodule of $M$. Note also that each highest weight module belongs to $\CO$, in particular, $\Delta(\lambda)$ and $L(\lambda)$ belong to $\CO$ for any $\lambda\in \hfhd$.

\subsection{Verma flags}

We say that an object $M$ of $\CO$ {\em admits a Verma flag} if there exists a finite filtration
$
0=M_0\subset M_1\subset \dots\subset M_n=M
$
such that for each $i=1,\dots, n$ the subquotient $M_{i}/M_{i-1}$ is isomorphic to $\Delta(\mu_i)$ for some $\mu\in \fhd$.
The number of occurences of $\Delta(\mu)$ in a filtration as above is in fact independent of the filtration and is denoted by $(M:\Delta(\mu))$.

\subsection{BGG-reciprocity}

Let $\CA$ be an abelian category and $L$ a simple object in $\CA$. Recall that a projective cover of $L$ in $\CA$ is a non-zero morphism $f\colon P\to L$ that has the property that $P$ is projective and that any morphism $g\colon Q\to P$ such that $f\circ g$ is non-zero, is surjective. Often one calls $P$ a projective cover of $L$ and assumes that the morphism $f$ is given.
The following result is due to Bernstein, Gelfand and Gelfand and is a characteristic zero version of a similar result of Humphreys.

\begin{theorem}[\cite{BGG}] Let $\lambda\in\fhd$.
\begin{enumerate}
\item 
There exists an up to isomorphism unique projective cover $P(\lambda)\to L(\lambda)$ in $\CO$.
\item 
The object $P(\lambda)$ admits a Verma flag and for the multiplicities holds the BGG-reciprocity formula
$
(P(\lambda):\Delta(\mu))=[\Delta(\mu):L(\lambda)].
$
\end{enumerate}
\end{theorem}

Using the BGG-reciprocity we can now reformulate the Kazhdan--Lusztig conjecture:

\begin{conjecture}\label{conj-KL3} Suppose that $\lambda$ is integral, regular and dominant. Then we have for all $y\in\CW$
$$
(P(y. \lambda):\Delta(\mu))=
\begin{cases} P_{w_0x,w_0y}(1), &\text{ if $\mu=x. \lambda$ for some $x\in\CW$},\\
0, & \text{ otherwise.}
\end{cases}
$$
\end{conjecture}

The reason for this additional reformulation of the Kazhdan--Lusztig conjecture is that both $P(\lambda)$ and $\Delta(\lambda)$ admit {\em deformed versions} and can hence be studied in a relative setting. This is not possible for the simple object $L(\lambda)$.
 
\subsection{Deformed category $\CO$}

The main idea of deformation theory is to not only consider a certain Verma module $\Delta(\lambda)$ by itself, but to study {\em families} of Verma modules.  Having understood all Verma modules in a neighbourhood of $\Delta(\lambda)$ we hope to obtain information on $\Delta(\lambda)$ itself. The method for this is as follows (references for the following are \cite{Soe2, FieMZ,FieCombofO}).

Let $S=S(\fh)$ be the symmetric algebra of the vector space $\fh$ and let $A$ be a commutative, unital, Noetherian, finitely generated $S$-algebra. We call such an algebra a {\em deformation algebra} in the following. We can now consider $\fg_A:=\fg\otimes_\DC A$ as an $A$-Lie algebra and study $\fg_A$-modules. Note that such an object is an $A$-module $M$ endowed with an $A$-linear action of $\fg$. 

As $A$ is supposed to be a unital $S$-algebra it  comes with the structure homomorphism $\tau\colon S\to A$, $f\mapsto f\cdot 1_A$. For $\lambda\in\fhd$ we define the $\fb_A$-module $A_\lambda$ as the free $A$-module of rank $1$ on which $H\in \fh$ acts as multiplication with the scalar $\lambda(H)+\tau(H)\in A$ and $[\fb,\fb]$ acts trivially. By induction we obtain the {\em deformed Verma module} 
$$
\Delta_A(\lambda):=U(\fg)\otimes_{U(\fb)} A_\lambda.
$$ 
The straightforward generalization of the definition of $\CO$ is the following:

\begin{definition} Let $M$ be a $\fg_A$-module. 
\begin{enumerate} 
\item $M$ is called a {\em weight module} if $M=\bigoplus_{\lambda\in\fhd} M_\lambda$, where $M_\lambda=\{m\in M\mid H.m=(\lambda(H)+\tau(H))m\text{ for all $H\in\fh$}\}$.
\item $M$ is called {\em locally $\fb_A$-finite}, if every $m\inÊM$ is contained in a $\fb_A$-submodule of $M$ that is finitely generated over $A$. 
\end{enumerate}
We denote by $\CO_A$ the full subcategory of the category of $\fg_A$-modules that consists of locally $\fb_A$-finite weight modules.
\end{definition}

\subsection{Simple objects in $\CO_A$}
From now on we assume that $A$ is a {\em local} deformation algebra with maximal ideal $\fm\subset A$ and quotient field $K=A/\fm$. Then we can consider $K$ as a deformation algebra as well. Note that $\fg_K$ is again a Lie algebra over a field of characteristic zero and all our "non-deformed" results apply. 

The corresponding category $\CO_K$ is just a direct summand of the usual category $\CO$ over the Lie algebra $\fg_K$ that contains all modules whose set of weights is contained in the image of the affine space $\tau+\fhd$ inside the $K$-linear dual of $\fh\otimes_\DC K$. In particular, each $\lambda\in \fhd$ parametrizes a simple object $L_K(\lambda)$ in $\CO_K$ with highest weight $\lambda+\tau$.

The locality of $A$ allows us to apply the Nakayama lemma for the proof of the following statement.

\begin{theorem} \cite{FieMZ} The base change functor $\cdot\otimes_A K\colon \CO_A\to \CO_K$ induces a bijection
$$
\left\{
\begin{matrix}
\text{simple isomorphism} \\
\text{classes of $\CO_A$}
\end{matrix}
\right\}\stackrel{\sim}\to
\left\{
\begin{matrix}
\text{simple isomorphism} \\
\text{classes of $\CO_K$}
\end{matrix}
\right\}.
$$
\end{theorem}

We now need to find deformed versions of the projective covers.

\subsection{Deformed projective objects}
As in the non-deformed case one proves the first part of the following theorem. The second part uses the idempotent lifting lemma.

\begin{theorem}\cite{FieMZ} Suppose that $A$ is a local deformation algebra with residue field $K$. Let $\lambda\in\fhd$.
\begin{enumerate}
\item 
There exists an up to isomorphism unique projective cover $P_A(\lambda)\to L_A(\lambda)$ in $\CO_A$.
\item 
The object $P_A(\lambda)$ admits a deformed Verma flag and for the multiplicities holds the BGG-reciprocity formula
$
(P_A(\lambda):\Delta_A(\mu))=[\Delta_K(\mu):L_K(\lambda)].
$
\end{enumerate}
\end{theorem} 

Note that on the right of the BGG-reciprocity above the Verma module and the simple object are defined over the residue field $K$. This allows us to give yet another formulation of the Kazhdan--Lusztig conjecture. Let $S=S(\fh)$ be the symmetric algebra of the vector space $\fh$ and denote by $\tS$ its localization at the maximal ideal $S\cdot \fh$. As the residue field of $\tS$ is $\DC$ we obtain the following equivalent formulation of the Kazhdan--Lusztig conjecture:

\begin{conjecture} \label{conj-projmul} Suppose that $\lambda\in\fhd$ is integral, regular and dominant. Then we have for all $w\in\CW$
$$
(P_{\tS}(w.\lambda):\Delta_{\tS}(\mu))=
\begin{cases}
P_{w_0x,w_0w}(1), &\text{ if $\mu=x. \lambda$ for some $x\in\CW$},\\
0, & \text{ otherwise.}
\end{cases}
$$
\end{conjecture}
Finally, we can relate the deformed category $\CO$ to the category of sheaves on a moment graph.

\subsection{A functor into moment graph combinatorics}

Let $w_0\in\CW$ be the longest element in the Weyl group. Then  {\em Soergel's structure functor} (\cite{Soe1}) is given by 
$$
\DV:=\Hom(P_A(w_0. \lambda), \cdot) \colon \CO_A\to \catmodrÊ\End(P_A(w_0. \lambda)).
$$

Let us consider the case $A=\tS$ again. Let $Q$ be the quotient field of $\tS$.  Each endomorphism $f$ of $P_{\tS}(w_0. \lambda)$ induces a homomorphism $f_Q$  of $P_{\tS}(w_0. \lambda)\otimes_{\tS} Q$. As  $P_{\tS}(w_0. \lambda)$ admits a deformed Verma flag, it is torsion free as a $\tS$-module, hence the map $f\mapsto f_Q$ is injective. So we can consider $\End_{\CO_{\tS}}(P_{\tS}(w_0. \lambda))$  as a submodule in $\End_{\CO_Q}(P_{\tS}(w_0. \lambda)\otimes_{\tS} Q)$.

\begin{proposition}[\cite{FieCombofO}]
\begin{enumerate}
\item  There is an isomorphism
$
 P_{\tS}(w_0. \lambda)\otimes_{\tS} Q\cong\bigoplus_{x\in\CW}\Delta_Q(x. \lambda).
$
It induces a {\em canonical} identification
$$
\End_{\CO_Q}(P_{\tS}(w_0. \lambda)\otimes_{\tS} Q)=\bigoplus_{x\in\CW} Q.
$$
\item 
The subspace $\End_{\CO_{\tS}}(P_{\tS}(w_0. \lambda))\subset \End_{\CO_Q}(P_{\tS}(w_0. \lambda)\otimes_{\tS} Q)=\bigoplus_{x\in\CW} Q$ coincides with
$$
\left\{
(z_x)\in \bigoplus_{x\inÊ\CW} \tSÊ\left|
\begin{matrix}
z_x\equiv z_{s_{\alpha} x} \mod \alpha^{\vee}\\
\text{ for all $x\in\CW$ and $\alpha\in R^+$}
\end{matrix}
\right.
\right\}.
$$
\end{enumerate}
\end{proposition}

Now let us consider $\DC$-sheaves on the moment graph associated to $R$. The definition of 
 the structure algebra $\CZ$ on this graph resembles greatly the result we obtained for $\End(P_{\tS}(w_0. \lambda))$, except for one  detail: the graded symmetric algebra is replaced by a localized symmetric algebra. We obtain the following result.

\begin{corollary} There is a {\em canonical} isomorphism $\End(P_{\tS}(w_0. \lambda))\cong \CZ\otimes_S \tilde S$.
\end{corollary}

We can now consider $\DV$ as a functor from $\CO_{\tS}$ to $\CZ\otimes_S \tilde S\catmod$. Let us denote by $\CO_{\tS,[\lambda]}$ the block of $\CO_{\tS}$ containing $L_{\tS}(\lambda)$ and by $\CO_{\tS,[\lambda]}^{Verma}$ its subcategory of objects admitting a Verma flag. It inherits an exact structure from the natural exact structure of the abelian category $\CO_{\tS}$. Here is now our link to moment graph theory:

\begin{theorem} [\cite{FieAdv}]\label{theorem-Vequiv}
\begin{enumerate}
\item The functor $\DV$ induces an {\em exact} equivalence between the  $\CO_{[\lambda]}^{Verma}$ and $\CC^{Verma}\otimes_S \tilde S$. 
\item For $w\in\CW$ there is an isomorphism 
$
\DV(P_{\tS}(w. \lambda))\cong \Gamma(\SB(w))\otimes_S \tilde S
$
of $\CZ\otimes_{S}\tS$-modules. 
\item We have for each $x\in\CW$
$
(P_{\tS}(w. \lambda):\Delta_{\tS}(x. \lambda))=\rk\, \SB(w)^x.
$
\end{enumerate}
\end{theorem}
Here, $\rk$ refers to the ungraded rank of the free $S$-module $\SB(w)^x$.
 Hence  Conjecture \ref{conj-projmul} is equivalent to 
\begin{conjecture} \label{conj-BMP} Suppose that $k=\DC$. Then we have for all $x,w\in\CW$
$$
\rk\,\SB(w)^x=P_{w_0x,w_0w}(1).
$$
\end{conjecture}
This (ungraded) conjecture is in fact equivalent to a graded version (cf.~\cite{FieTAMS}).

\section{Geometry}\label{sec-top}

Now we want to relate the geometry of flag varieties to the Braden-MacPherson sheaves on the underlying Bruhat graph. A reference for the following is  \cite{FW}, which contains a positive characteristic analog of the main result in \cite{GKM}. For the definition of the equivariant derived category, see \cite{BL}.

\subsection{$H$-spaces}

Let $H$ be a topological group. An {\em $H$-space} is a topological space $X$ together with a continuous $H$-action $H\times X\to X$. An $H$-space is called {\em (topologically) free} if the quotient map $X\to X/H$ is an $H$-bundle (i.e.~ a locally trivial fibration with fiber $H$). 

For the following we fix an $H$-space $EH$ that is free and contractible (such a space always exists).  This space is not uniquely defined, but nothing that follows depends on our choice for $EH$. Then we can define for an arbitrary $H$-space  the topological space 
$$
X_H:=X\times_H EH,
$$ 
which is the orbit space of $X\times EH$ under the diagonal $H$-action. We then have a diagram

\centerline{
\xymatrix{
& X\times EH \ar[ld]_q\ar[rd]^p&\\
X_H&&X,
}
}
 \noindent
 where $q$ is the canonical orbit map and $p$ is the projection onto the first factor.

\subsection{The equivariant derived category}
For the sheaf theory we now also need a field $k$ of coefficients. To the data $(H,X,k)$ one associates the following category:

\begin{definition} The {\em equivariant derived category} $D^+_H(X,k)$ of sheaves on $X$ with coefficients in $k$ is the full subcategory  of $D^+(X_H,k)$ that contains all sheaves $\CF$ for which there is a sheaf $\CF_X\in D^+(X,k)$ such that $q^\ast(\CF)\cong p^\ast(\CF_X)$.
\end{definition}
By $D^+(Y,k)$ we denote the derived category of sheaves of $k$-vector spaces on a topological space $Y$ with cohomology bounded from below.

Suppose that $X$ and $X^\prime$ are $H$-spaces and that $f\colon X\to X^\prime$ is an $H$-equivariant map. Then $f\times \id\colon X\times EH\to Y\times EH$ induces a continuous map $f_H\colon X_H\to Y_H$ and (for suitable $f$, see \cite{BL}) we get base change functors $f_H^\ast$, $f_H^!$, $f_{H\ast}$ and $f_{H!}$ between $D^+(X_H,k)$ and $D^+(Y_H,k)$. One checks, again under suitable assumptions on the map $f$, that these functors induce functors between the subcategories $D^+_H(X,k)\subset D^+(X_H,k)$ and $D^+_H(Y,k)\subset D^+(Y_H,k)$. We denote these restrictions by $f^\ast$, $f^!$, $f_\ast$ and $f_!$ in order to save indizes.

\subsection{Hypercohomology}

The map $\pi\colon X\to \{pt\}$ induces a direct image functor
$
\pi_{\ast}\colon D^+_H(X,k)\to D^+_H(\pt,k).
$
The category $D^+_H(\pt,k)$ is a full subcategory of $D^+(BH,k)$, where $BH=EH/H=\pt\times_H EH$ is the classifying space of $H$. We denote by  
$$
A_H:=H^\ast(BH,k)
$$ 
its ordinary cohomology with coefficients in $k$.

\begin{definition} Let $\CF\in D^+_H(X,k)$. The {\em equivariant hypercohomology} of $\CF$ is the graded $A_H$-module
$
\HHypgr(\CF):=H^\ast(\pi_\ast \CF),
$, i.e.~  the (ordinary) cohomology of the space $BH$ with coefficients in the sheaf $\pi_\ast\CF$.
\end{definition}

Suppose that $i\colon Y\to X$ is the inclusion of a $H$-stable subvariety.  For any sheaf $\CF\in D^+_H(Y,k)$ we denote by
 $
 \CF_Y:= i^\astÊ\CF
 $
 its restriction to $Y$. Note that this is an equivariant sheaf on $Y$. The adjunction $\id\to i_\ast i^\ast$ yields a morphism 
 $
 \CF\to i_\ast i^\astÊ\CF=i_\ast \CF_Y
$
between sheaves on $X$. Applying the hypercohomology functor $\HHypgr$ yields the {\em restriction morphism}
$$
\HHypgr(\CF)\to \HHypgr(\CF_Y)
$$
between $A_H$-modules.

\subsection{The main example - flag varieties as $T$-spaces}

Let $G$ be a connected reductive complex algebraic group with root system $R$, and let $T\subset B\subset G$ be a maximal torus inside a Borel subgroup of $G$. The quotient $X=G/B$ carries a canonical structure of a (projective) algebraic variety. It is acted upon (algebraically) by the torus $T$. 

Note that each complex algebraic variety can be viewed as a topological space with the underlying {\em metric} topology. In particular, we can view $T$ as a topological group and $G/B$ as a $T$-space. 

Choose an isomorphism $T\cong (\DC^\times)^r$. Then $T$ acts by componentwise multiplication on the space $(\DC^n\setminus\{0\})^r$. We embed $\DC^n\setminus\{0\}$ into $\DC^{n+1}\setminus\{0\}$ by adding a $0$ on the $(n+1)^{\text{st}}$ coordinate. The space $(\DC^\infty\setminus \{0\})^r:=\lim_{n\to\infty} (\DC^n\setminus\{0\})^r$ is then  a contractible space with a topologically free $T$-action, hence we can take this as a model for $ET$. 

 Note that $BT=(\DC^\infty\setminus\{0\})^r/T=\DP^{\infty}$, so $A_T$ can be identified with the symmetric algebra $S=S(X^\ast(T) \otimes_\DZ k)$ (here $X^\ast(T)=\Hom(T,\DC^\times)$ is the character lattice of $T$).

\subsection{The moment graph associated to a $T$-space $X$}

From now on we restrict ourselves to the case that the group $H$ is a complex torus $T$.
One of many possible methods to calculate the hypercohomology $\THypgr(\CF)$ of 
a $T$-equivariant sheaf $\CF$ is the localization method of Goresky, Kottwitz and MacPherson. In the following we will shortly review their method.   

Each $\alpha\in X^{\ast}(T)$ defines an action of $T$ on the affine variety $\DC^\times$ by $t.x=\alpha(t)x$ for $t\in T$, $x\in \DC^\times$. We denote the resulting $T$-space by $\DC^\times_\alpha$. 
In order to associate a moment graph to a $T$-variety $X$ we now assume the following.

\begin{enumerate}
\item There are only finitely many $0$- and $1$-dimensional $T$-orbits in $X$.
\item Each fixed point is attractive (recall that  a fixed point $x
\in X$ is called attractive if all weights of $T$ on the tangent space
$T_xX$ lie in an open half space of $X^\ast\otimes_\DZ \DR$). 
\item The closure of a $1$-dimensional orbit in $X$ is smooth.
\item For each $1$-dimensional orbit $E$ in $X$ there is a character $\alpha\in X^{\ast}(T)$ such that $E\cong \DC^\times_{\alpha}$ as a $T$-space.
\end{enumerate}
Note that the character $\alpha_E$ in part (3) is only well-defined up to a sign. For each $1$-dimensional orbit $E$ we now fix $\alpha_E$. Nothing that follows depends on this choice. 

The moment graph $\CG_X$ associated to a $T$-variety $X$ that satisfies the assumptions above is the following. Its set of vertices is the set $X^T$ of fixed points of $X$. The vertices $x,y\in X^T$, $x\ne y$ are connected by an edge if there is a $1$-dimensional orbit $E\subset X$ such that $\ol E=E\cup\{x,y\}$. We denote this edge by $E$ as well and we set $\alpha(E):=\alpha_E$. So $\CG_X$ is a moment graph over the weight lattice $X^\ast(T)$.
 
 If $X=G/B$ is the flag variety associated with the root system $R$, then the moment graph $\CG_X$ is the following: The set of fixed points can be identified with the set $\CW$, and $x$ and $y$ are connected by an edge if $x=s_\alpha y$ for some $\alpha\in R$. We label this edge by $\alpha$. So in contrast to the example in  Section \ref{subsec-momgra}, the labels are given by roots and not by coroots! This is the reason why our localization relates the representation theory of $\fg$  to the topology of the Langlands dual flag variety $G^\vee/B^\vee$ rather than $G/B$.
 
 \subsection{Moment graph sheaves associated to equivariant sheaves}
 Our next step is to associate a $k$-sheaf $\DW(\CF)$ on $\CG_X$ to any $\CF\in D_T^+(X,k)$. For a vertex $x\in X^T$ we set
 $$
 \DW(\CF)^x := \THypgr(\CF_x),
 $$
 and for a one dimensional orbit $E$ we set 
 $$
 \DW(\CF)^E:= \THypgr(\CF_E).
 $$
From the identification $E\cong \DC^\times_\alpha$ one deduces that $\alpha\THypgr(\CF_E)=\{0\}$, as required.
 For the construction of the homomorphisms $\rho_{x,E}$ we need a little lemma.

Let $Z$ be a  $T$-variety with an attractive fixed point $x$. 
 We denote by $i\colon \{x\}\to Z$ the inclusion and by $\pi\colon Z\to\{x\}$ the projection. If we apply the functor $\pi_\ast$ to the adjunction $\id\to i_{\ast}i^\ast$, we get a natural morphism
$
\pi_{\ast}\to i^\ast
$
(since $\pi\circ i$ is the identity map). 

\begin{lemma}[\cite{FW}] Suppose that $Z$ is connected and affine and let $x\in Z$ be an attractive fixed point. Then the morphism $\pi_{\ast}\to i^\ast$ constructed above is an isomorphism.
\end{lemma}

Now let us consider again the case of a non-necessarily affine $T$-space $X$ as before. Let $\CF$ be an equivariant sheaf on $X$, let $E\subset X$ be a  one dimensional orbit and choose a fixed point $x$ in the closure of $E$. From the above lemma we obtain an isomorphism
$
\pi_{\ast}\CF_{E\cup x}\stackrel{\sim}\to i^\ast\CF_{E\cup x}=\CF_x
$
of sheaves on the point $\{x\}$, and after taking hypercohomology we get an isomorphism 
$
\THypgr(\CF_{E\cup x})\stackrel{\sim}\to \THypgr(\CF_{x}).
$
We can  now define $\rho_{x,E}\colon \THypgr(\CF_{x})\to \THypgr(\CF_E)$ as the composition of the invers of the above morphism with the restriction morphism:
$$
\rho_{x,E}\colon \THypgr(\CF_{x}) \stackrel{\sim}\leftarrow
\THypgr(\CF_{E\cup x})\to \THypgr(\CF_{E}).
$$

Hence we constructed the remaining ingredient for a moment graph sheaf. The constructions above are clearly functorial,  so we now have a functor
$$
\DW\colon D_T^+(X,k)\to \CG_X\catmod_k.
$$

\subsection{The localization theorem}
The next result shows that one can recover the  global hypercohomology $\THypgr(\CF)$ of certain equivariant sheaves $\CF$ from the local hypercohomologies on fixed points and one dimensional orbits. In the case $\ch\, k=0$, this is  a result of Goresky, Kottwitz and MacPherson (\cite{GKM}). With some additional care one can use their arguments in order to prove the statement  for almost arbitrary characteristic (\cite{FW}). 

\begin{theorem}[\cite{GKM,FW}] Suppose that $(\CG_X,k)$ satisfies the GKM-assumption and suppose that $\CF\in D^+_T(X,k)$ is such that $\THypgr(\CF)$ is a free $S$-module. Then we have an isomorphism
$$
\THypgr(\CF)=\Gamma(\DW(\CF)).
$$
\end{theorem}

Of course, now we should be looking for  equivariant sheaves on $X$ that correspond to the Braden-MacPherson sheaves on $\CG$. The answer may come as a surprise: these are not the intersection cohomology sheaves, but the {\em parity sheaves}.

\subsection{(Equivariant) parity sheaves on stratified varieties}

For the definition of parity sheaves we need yet another piece of data, namely a {\em stratification} of the variety $X$. Recall that a stratification of $X$ is a decomposition 
$$
X=\bigsqcup_{\lambda\in\Lambda} X_\lambda
$$
by locally closed subvarieties $X_\lambda\subset X$ such that the closure of each stratum is  a union of strata.  For any $\lambda\in\Lambda$ we denote by $i_\lambda\colon X_\lambda\to X$ the inclusion of the stratum $X_\lambda$. We furthermore impose the following assumptions:

\begin{itemize}
\item Each stratum is  $T$-stable and there is a $T$-equivariant isomorphism $X_\lambda\to \DC^{n_\lambda}$, where $\DC^{n_\lambda}$ carries a linear $T$-action.
\item There are only finitely many $0$- and $1$-dimensional orbits in $X$.
\item The first two assumptions imply that there is a unique fixpoint $x_\lambda$ in $X_\lambda$ for all $\lambda$. We assume that $x_\lambda$ is attractive. 
\item We assume that the stratification is a Whitney-stratification.
\end{itemize}

In the case of our main example, i.e.~  the flag manifold $X=G/B$, we consider  the stratification given by $B$-orbits. By the Bruhat decomposition we can identify the set of orbits with the Weyl group $\CW$.

\begin{definition}[\cite{JMW}] Let $\CP\in D^+_T(X,k)$.
\begin{itemize}
\item $\CP$  is called {\em even}, if for all $\lambda\in \Lambda$ the sheaves $i_\lambda^\ast\CP$ and $i_\lambda^!\CP$ are isomorphic to a direct sum of constant sheaves  shifted by even degrees.
\item $\CP$ is called   {\em odd}, if $\CP[1]$ is even.
\item $\CP$ is called {\em parity}, if it is a direct sum of an even and an odd sheaf.
\end{itemize}
\end{definition}

Parity sheaves do not always exist. In the case of flag varieties we have the following result.

\begin{proposition}[\cite{JMW}] For each $w\in\CW$ there exists an up to isomorphism unique indecomposable parity sheaf $\CP(w)$ on the flag variety $G/B$ with $\supp\,\CP(w)\subset \ol{X_w}$ and $\CP(w)_w\cong\ul{k}_{X_w}$.
\end{proposition}

Now we can link the geometry of flag varieties to the Braden-MacPherson sheaves defined before:

\begin{theorem}[\cite{FW}] For each $w\in \CW$ there is an isomorphism 
$$
\DW(\CP(w))\cong \SB(w).
$$
In particular, for each pair $x,w\in\CW$ we have an isomorphism
$
\THypgr(\CP(w)_x)\cong\SB(w)^x.
$
\end{theorem}

We need one last step in order to prove the Kazhdan--Lusztig conjecture. The {\em decomposition theorem} of Beilinson, Bernstein, Deligne and Gabber (cf.~\cite{BBD}) is used to prove the following:
\begin{theorem} Suppose that $\ch k=0$. Then $\CP_k(w)$ is the intersection cohomology complex $\IC(\ol{BwB/B},k)$ on the Schubert variety $\ol{BwB/B}$.
\end{theorem}

Kazhdan and Lusztig managed to calculate the ranks of the stalks of the intersection cohomology complexes:

\begin{theorem}[\cite{KL1,KL2}]\label{theorem-multIC} Let $x,w\in\CW$. Then 
$$
\rk\, \THypgr(IC(\ol{BwB/B}),k)_x)=P_{w_0x,w_0w}(1).
$$
\end{theorem}

In particular, we obtain $\rk \SB(w)^x=h_{x,w}(1)$ for the Braden--MacPherson sheaves on $\CG_{G/B}$. The Kazhdan--Lusztig polynomials depend only on the underlying Coxeter system, not on the root system, so they do not change when we switch to the Langlands dual setup. So we also obtain a proof of Conjecture \ref{conj-BMP} and hence of the Kazhdan--Lusztig Conjecture.

\subsection{The Elias--Williamson work}
Very recently, the remarkable paper \cite{EW} appeared. It contains a direct, i.e.~ non-geometric proof of Conjecture \ref{conj-BMP} in characteristic $0$ and hence of the Kazhdan--Lusztig conjecture. The proof uses Soergel bimodules instead of moment graph sheaves. But by a result in \cite{FieTAMS}, the category of Soergel bimodules is equivalent to the category of Braden--MacPherson sheaves, and Conjecture \ref{conj-BMP} translates into which is known as Soergel's conjecture. Unfortunately, the arguments used by Elias and Williamson do not generalize to the positive characteristic case. 

\end{document}